# A CONCISE AND DIRECT PROOF OF "FERMAT'S LAST THEOREM"

## by Roger Ellman

*FERMAT'S LAST THEOREM* states:

There can be no non-zero integer solution for $n>2$ to the equation

$(1) \quad a^n + b^n = c^n$

### Step 1

Restate the problem as follows:

For $x$, $i$, $n$ and $f(x,i)$ all non-zero integers and $i<x$ there is no solution for $n>2$ to the equation

$(2) \quad x^n = [x-i]^n + [f(x,i)]^n$

That is, make the following substitutions in equation (1):

$$x^n = c^n \qquad [x-i]^n = a^n \qquad [f(x,i)]^n = b^n$$

Clearly there is no difficulty with the $x^n$ term nor the $[x-i]^n$ term. Both are integers and perfect $n^{th}$ powers of integers. The issue now is:

Can $f(x,i)$ be a non-zero integer for $n>2$ and equation (2) still valid?

### Step 2

The 1st constraint on $b^n$: it must be the difference of $c^n$ and $a^n$.

$(3) \quad b^n = [f(x,i)]^n$

$\qquad = x^n - [x-i]^n \qquad$ [Solving equation (2)]

$\qquad = x^n - [x^n - n \cdot x^{n-1} \cdot i + \ldots \pm i^n] \qquad$ [Binomial expansion]

$\qquad = n \cdot x^{n-1} \cdot i - \ldots \pm i^n$

### Step 3

The 2nd constraint on $b^n$: it must be a perfect $n^{th}$ power.

$(4) \quad b^n = [x-j]^n = [x-j]_1 \cdot [x-j]_2 \cdot [x-j]_3 \cdot \cdots \cdot [x-j]_n$

$\qquad$ where: $b = x-j$ (just as $a = x-i$)
$\qquad \qquad \qquad j$ is a non-zero integer, $j<x$

### Step 4

These two constraints are simultaneous. They are for the same $b^n$. Therefore the two expressions must be identical; they must always simultaneously deliver the same value of $b^n$.

The order of Step 2, equation (3) is one less than the order of Step 3, equation (4). To compare the two expressions as an identity their order must be the same. That is accomplished by removing one factor of $b$ from each of equations (3) and (4), as follows.

$(5)\quad b^n = n \cdot x^{n-1} \cdot i - \ldots \pm i^n$ [equation 3]

$$= \underbrace{\frac{n \cdot i}{m}}_{b} \cdot \underbrace{m \cdot \left[ x^{n-1} - \ldots \pm \frac{i^{n-1}}{n} \right]}_{b^{n-1}}$$

(The parameter $m$ is necessary because the quantity, $n \cdot i$, which factored out normalizes the expressing, is not necessarily equal to $b$.)

$(6)\quad b^n = \underbrace{[x-j]_1}_{b} \cdot \underbrace{[x-j]_2 \cdots [x-j]_n}_{b^{n-1}}$ [equation 4]

$$= \underbrace{[x-j]_1}_{b} \cdot \underbrace{m \cdot \left[ [x-k]_2 \cdot [x-k]_3 \cdot \cdots \cdot [x-k]_n \right]}_{b^{n-1}}$$

(The $m$ here is for identity to be possible -- for the coefficient of the $x^{n-1}$ term in the two expressions to be able to be equal, when $m \neq 1$.)

Step 5

Now, expression (5) and expression (6) must yield the same value for $b^n$ for all values of $x$. To establish that condition we will require, for convenience rather than the entire expressions, that $[b^{n-1}/m]$ in each expression yield the same value for all values of $x$. The two expressions are (using the binomial theorem expansion formula): in expression (5)

$(7)\quad x^{n-1} - \frac{[n-1]}{2 \cdot 1} \cdot x^{n-2} i + \frac{[n-1][n-2]}{3 \cdot 2 \cdot 1} \cdot x^{n-3} i^2 - \ldots \pm \frac{i^{n-1}}{n}$

and in expression (6)

$(8)\quad x^{n-1} - \frac{[n-1]}{+1} \cdot x^{n-2} k^1 + \frac{[n-1][n-2]}{2 \cdot 1} \cdot x^{n-3} k^2 - \ldots \pm k^{n-1}$

Equating the pair of terms of zero order in equations (7) and (8):

$(9)\quad \pm \frac{i^{n-1}}{n} = \pm k^{n-1}$

$$k = \frac{i}{\sqrt[n-1]{n}}$$

The $[n-1]^{th}$ root of $n$ is irrational for $n>2$. Therefore, for $n>2$, $k$ is irrational and $b$ is irrational and cannot be an integer, which proves the theorem.

However, $k$ in expression (8) is a function of $x$. The only values of $k$ that are able to make the expression for $b^{n-1}$ in the horizontal bracket to the right in the second line of expression (6) actually be equal to $b^{n-1}$ are as follows:

(10)
$$k = \left[x - \left[\frac{b^n}{n \cdot i}\right]^{1/[n-1]}\right] \qquad \text{[where } b \text{ is also a function of } x\text{]}$$

which can readily be verified by substitution, that is

(11)
$$m \cdot [x-k]^{n-1} = \frac{n \cdot i}{b} \cdot \left[x - \overbrace{\left[x - \left[\frac{b^n}{n \cdot i}\right]^{1/[n-1]}\right]}^{k}\right]^{n-1}$$

$$= \frac{n \cdot i}{b} \cdot \left[\left[\frac{b^n}{n \cdot i}\right]^{1/[n-1]}\right]^{n-1} = \frac{n \cdot i}{b} \cdot \frac{b^n}{n \cdot i} = b^{n-1}$$

The problem with $k$ being a function of $x$ is that the apparent terms of given orders of $x$ and their coefficients are not necessarily as they appear in expression (8) when expression (9) is substituted for $k$ in expression (8). However, if the term coefficients experience no net change from the substitution, then the comparison of any pair of coefficients is valid even though $k = f(x)$. That is exactly the situation in the present case (and may relate to why the theorem withstood proof for three centuries) as follows.

To show this in an overall general form would be too algebraically complex to contemplate. The pattern can be developed with two examples.

<u>Example #1:   n = 2</u>

| Expression Nr As on Page 2 | Content |
|---|---|
| (5) | $b^n = 2 \cdot x \cdot i - i^2$ |
|     | $= \frac{2 \cdot i}{m} \cdot m \cdot \left[x - \frac{i}{2}\right]$ |
| (6) | $b^n = [x-j] \cdot [x-j]$ |
|     | $= [x-j] \cdot m \cdot [x-k]$ |
| (7) | $[b^{n-1}/m] = x - i/2$ |
| (8) | $[b^{n-1}/m] = x - k$ |
| (10) Page 3 | $k = \left[x - \left[\frac{b^2}{2 \cdot i}\right]^{1/1}\right]$ |
|     | $= \left[x - \left[\frac{2 \cdot x \cdot i - i^2}{2 \cdot i}\right]^{1/1}\right]$ |
|     | $= i/2$ |
| Substituting (10) For the $k$ in (8) gives (8) ≡ (7) | $\left[b^{n-1}/m\right] = x - i/2$ |

|                              | Example #2:  n = 3                                              |
| ---------------------------- | --------------------------------------------------------------- |
| Expression Nr As on Page 2   | Content                                                         |

(5)  $b^n = 3 \cdot x^2 \cdot i - 3 \cdot x \cdot i^2 + i^3$

$\phantom{b^n} = \dfrac{3 \cdot i}{m} \cdot m \cdot \left[ x^2 - x \cdot i + \dfrac{i^2}{3} \right]$

(6)  $b^n = [x-j] \cdot [x-j] \cdot [x-j]$

$\phantom{b^n} = [x-j] \cdot m \cdot [x-k] \cdot [x-k]$

(7)  $[b^{n-1}/m] = x^2 - x \cdot i + i^2/3$

(8)  $[b^{n-1}/m] = x^2 - 2 \cdot k \cdot x + k^2$

(10) Page 3  $k = \left[ x - \left[ \dfrac{b^3}{3 \cdot i} \right]^{1/2} \right]$

$\phantom{k} = \left[ x - \left[ \dfrac{3 \cdot x^2 \cdot i - 3 \cdot x \cdot i^2 + i^3}{3 \cdot i} \right]^{1/2} \right]$

$\phantom{k} = x - [x^2 - x \cdot i + i^2/3]^{1/2}$

Substituting (10) For the $k$ in (8) gives (8) ≡ (7)

$\left[ b^{n-1}/m \right] = x^2 - x \cdot i + i^2/3$

This pattern persists for all positive integer values of $n$. Therefore, the term coefficients experience no net change from the substitution and the comparison of any pair of coefficients is valid even though $k = f(x)$. Therefore, expression (9) is valid and expression (9) shows that $k$, and therefore $b$, are irrational for $n>2$, which proves the theorem.